\documentclass[a4paper,11pt,twoside,leqno]{article}

\usepackage{amsmath,amsthm,amsfonts,latexsym,amscd,amssymb}
\numberwithin{equation}{section}
\input xypic
\def\qed{{\hbadness=10000\hfill\ \vbox{\hrule height.09ex
   \hbox{\vrule width.09ex height1.55ex depth.2ex \kern1.8ex
   \vrule width.09ex height1.55ex depth.2ex}\hrule height.09ex}\break
   \bigskip}}
\newtheorem{theorem}{Theorem}[section]

\theoremstyle{definition}

\theoremstyle{remark}

\newcommand{\n}{\noindent}

\begin{document}

\linespread{1}\title {\textbf{Some characterizations of rectifying curves on a smooth surface in Euclidean 3-space}}

\author{Akhilesh Yadav\thanks{Corresponding author}, Buddhadev Pal}
\date{}
\maketitle

\noindent\textbf{Abstract:} In this paper, we investigate sufficient condition for the invariance of a rectifying curve on a smooth surface immersed in Euclidean 3-space under isometry by using Darboux frame $\left\lbrace T, P, U\right\rbrace$. Further, we find the deviations of the position vector of a rectifying curve on the smooth surface along any tangent vector $T = a\phi_u + b\phi_v$ with respect to the isometry. We also find the deviations of the position vector of a rectifying curve on the smooth surface along the unit normal $U$ to the surface and along $P (= U \times T)$ with respect to the isometry.

\n\textbf {Mathematics Subject Classification 2020:} 53A04, 53A05, 53A15.

\n\textbf{Key words:} Isometry, Frenet-frame, Darboux-frame, rectifying curve.

\section{Introduction}
 \n In the study of space curve in the Euclidean 3-space $E^3$ at every point of the curve one can associate the moving orthonormal frame called Serret-Frenet frame $\left\lbrace T, N, B\right\rbrace$, consisting of the unit tangent vector, the principal normal vector and the binormal vector, respectively. According to Serret-Frenet frame at  every point of the curve there exist three orthogonal planes, so called rectifying, normal and osculating planes. Rectifying curves are introduced by B. Y. Chen in [1] as a space curve whose position vector always lies in its rectifying plane. Here, the rectifying plane is spanned by the tangent vector $T(s)$ and the binormal vector $B(s)$. Thus, the position vector $\gamma(s)$ of a rectifying curve $\gamma$ in $E^3$ satisfies the equation  $\gamma(s) = \lambda(s)T(s) + \mu(s)B(s)$ for some differentiable functions $\lambda(s)$ and $\mu(s)$. The rectifying curves are also studied in [3] as the extremal curves. 
 Further, authors studied rectifying curves via the dilation of unit speed curves on the unit sphere $S^2$ in the Euclidean 3-space and obtained a necessary and sufficient condition for which the centrode $d(s)$ of a unit speed curve $\gamma(s)$ is a rectifying curve [5]. Also in [7], authors defined a rectifying curve in the Euclidean $4$-space as a curve whose position vector always lies in orthogonal complement $N^\bot$ of its principal normal vector field $N$. \\
 
 \n In [10], authors studied rectifying curve on a smooth surface and obtained a sufficient condition for which a rectifying curve on a smooth surface remains invariant under isometry of surfaces by using Serret-Frenet frame. On the other hand, when we study space curve on a smooth surface immersed in Euclidean 3-space at every point of the curve another moving orthonormal frame called Darboux frame $\left\lbrace T, P, U\right\rbrace$ comes naturally, where $T$ is the unit tangent vector to the curve at that point, $U$ being the unit normal to the surface and $P = U\times T$. In [4], authors gave some characterizations of position vector of a unit speed curve in a regular surface immersed in Euclidean 3-space which always lies in the planes spanned by $\left\lbrace T, U\right\rbrace$, $\left\lbrace T, P\right\rbrace$ and $\left\lbrace P, U\right\rbrace$, respectively by using the Darboux frame. Thus motivated sufficiently, we study rectifying curve on a smooth surface immersed in Euclidean 3-space and investigate the sufficient condition for the invariance of a rectifying curve on the smooth surface under isometry by using Darboux frame instead of Frenet frame. The paper is arrange as follws: In section 2, we discuss some basic theory of unit speed parametrized curve on a smooth surface. Section 3 is devoted to the investigation of the sufficient condition for the invariance of a rectifying curve on a smooth surface immersed in Euclidean 3-space under isometry. In this section, we also find the deviations of the position vector of a rectifying curve on the smooth surface along any tangent vector $T = a\phi_u + b\phi_v$, the unit normal $U$ to the surface and along $U \times T$ with respect to the isometry.
 
\section{Preliminaries}
\n Let $\gamma: I \rightarrow E^3$, where $I = (\alpha, \beta)\subset\mathbb{R}$, be the unit speed parametrized curve that has at least four continuous derivatives. Then the tangent vector of the curve $\gamma$ be denoted by $T$ and given by $T(s) = \gamma^{'}(s)$, $\forall s \in I$,  where $\gamma^{'}$ denote the derivative of $ \gamma$ with respect to the arc length
parameter $s$. The binormal vector $B$ is defined by $B = T\times N$, where $N$ is the principal normal vector to the curve $\gamma$. The Frenet-Serret equations are given by
\begin{equation}
T^{'}(s) = \kappa(s)N(s),
\end{equation}
\begin{equation}
N^{'}(s) = -\kappa(s)T(s) + \tau(s)B(s),
\end{equation}
\begin{equation}
B^{'}(s) = -\tau(s)N(s),
\end{equation}
where $\kappa(s)$ and $\tau(s)$ are smooth functions of $s$, called curvature and torsion of the curve $\gamma$.

\n Let $\phi: V\subset\mathbb{R}^2\rightarrow S$ be the coordinate chart for a smooth surface $S$ immersed in Euclidean space $E^3$ and the unit speed parametrized curve $\gamma: I \rightarrow S\subset E^3$, where $I = (\alpha, \beta)\subset\mathbb{R}$, contained in the image of a surface patch $\phi$ in the atlas of $S$. Then $\gamma(s)$ is given by 
\begin{equation}
\gamma(s) = \phi(u(s), v(s)),\indent \forall s \in I.
\end{equation}
Now, the curve $\gamma(s)$ lies on the surface $S$ there exists another moving orthonormal frame called Darboux frame $\left\lbrace T, P, U\right\rbrace$ at each point of the curve $\gamma(s)$. Since the unit tangent $T$ is common in both Frenet frame and Darboux frame, the vectors $N, B, P, U$ lie in the same plane. So that the relations between these frames can be given as follows:
 \begin{equation}
\left[ {\begin{array}{c}
	 T \\
	 P \\
	 U \\
	\end{array} } \right]=
 \left[ {\begin{array}{ccc}
 	1 & 0 & 0 \\
 	0 & cos\theta & sin\theta\\
 	0 & -sin\theta & cos\theta\\
 	\end{array} } \right]
 \left[ {\begin{array}{c}
 	T \\
 	N \\
 	B \\
 	\end{array} } \right],
 \end{equation}
 where $\theta$ is the angle between vectors $N$ and $P$.
 
 \n Again, since $\gamma(s)$ is unit-speed curve lies on surface $S$, $\gamma^{''}$ is perpendicular to $\gamma^{'}(= T)$, and hence is a linear combination of $U$ and $P (= U \times T)$. Thus 
 \begin{equation}
 \gamma^{''}(s) = k_n(s) U(s) + k_g(s) P(s),
 \end{equation}
 where $k_n$ and $k_g$ are smooth functions of $s$, called the normal curvature and the geodesic  curvature of $\gamma$, respectively. Since $U$ and $P$ are perpendicular unit vectors therefore from (2.6), we get
 \begin{equation}
 k_n(s) = \gamma^{''}(s).U(s)\indent and \indent k_g(s) = \gamma^{''}(s).P(s).
 \end{equation}
 Also from (2.1) and (2.7), we obtain
  \begin{equation}
 k_n(s) = \kappa(s)N(s).U(s)\indent and \indent k_g(s) = \kappa(s)N(s).P(s),
 \end{equation}
 which implies
 \begin{equation}
 k_n(s) = \kappa(s)sin\theta\indent and \indent k_g(s) = \kappa(s)cos\theta.
 \end{equation}
 Thus the curve $\gamma$ is a geodesic curve if and only if $k_g = 0$ and the curve $\gamma$ is an asymptotic line if and only if $k_n = 0$.
 
 \n Now, Differentiating (2.4) with respect to $s$, we get
 \begin{equation}
 T(s)= \gamma^{'}(s) = u^{'}\phi_u + v^{'}\phi_v.
 \end{equation}
 The unit normal $U$ to the surface $S$ is given by
 \begin{equation}
 U(s) = \frac{\phi_u\times\phi_v}{||\phi_u\times\phi_v||} = \frac{\phi_u\times\phi_v}{\sqrt{EG-F^2}}.
 \end{equation}
 Also, since $P = U \times T$ by using (2.10) and (2.11), we obtain
  \begin{equation}
 P(s) = \frac{1}{\sqrt{EG-F^2}}(Eu^{'}\phi_v + F(v^{'}\phi_v - u^{'}\phi_u) -G v^{'}\phi_u),
 \end{equation}
 where $E = \phi_u.\phi_u$, $F = \phi_u.\phi_v$ and $G = \phi_v.\phi_v$ are coefficients of first fundamental form.
 
\n\textbf{Definition 2.1.} [6] A diffeomorphism $f: S\rightarrow \overline{S}$ between two surfaces $S$ and $\overline{S}$ is an isometry if $<f_*(\omega_1), f_*(\omega_2)>_{f(p)} = <\omega_1, \omega_2>_{p}$, for all $p\in S$ and for all $\omega_1,\omega_2 \in T_pS$. The surfaces $S$ and $\overline{S}$ are called isometric if there is an isometry between them.

\n Now, from Theorem 5.1 and Corollary 8.2 of [9], we have the following:
 
 (i) An isometry $f$ between surfaces $S$ and $\overline S$ takes the geodesies of one surface to the geodesies of the other,
 
 (ii) Coefficients of first fundamental form preserve under isometry between surfaces $S$ and $\overline S$, i.e. if $E$, $F$, $G$ and  $\overline E$,   $\overline F$,  $\overline G$ are coefficients of first fundamental form of surfaces $S$ and $\overline S$, respectively then
\begin{equation}
E = \overline E,\indent F = \overline F\indent  and \indent G = \overline G.
\end{equation}

\section{Rectifying curves according to Darboux frame}

\n In this section, we study rectifying curves on a smooth surface by using Darboux frame. A curve $\gamma(s)$  on a smooth surface $S$ immersed in Euclidean 3-space  is called rectifying curve if its position vector always lies in rectifying plane of the curve. Thus the position vector of the curve $\gamma$ satisfies the equation 
\begin{equation}
\gamma(s) = \lambda(s)T(s) + \mu(s)B(s),
\end{equation}
for some differentiable functions $\lambda(s)$ and $\mu(s)$. Thus by using (2.5) in (3.1), we obtain
\begin{equation}
\gamma(s) = \lambda(s)T(s) + \mu(s)P(s)sin\theta + \mu(s)U(s)cos\theta.
\end{equation}
\n Now, from (2.9), (2.10), (2.11), (2.12) and (3.2), we get
\begin{equation}\begin{split}
\gamma(s)& = \lambda(s)(u^{'}\phi_u + v^{'}\phi_v) + \frac{\mu(s)k_g(s)}{k(s)\sqrt{EG-F^2}}(\phi_u\times\phi_v) \\& + \frac{\mu(s)k
_n(s)}{k(s)\sqrt{EG-F^2}}(Eu^{'}\phi_v + F(v^{'}\phi_v - u^{'}\phi_u) -G v^{'}\phi_u). 
\end{split}\end{equation}
This equation of rectifying curve on a smooth surface, which is neither a geodesic curve nor an asymptotic line on the surface. 

\n Now, if the rectifying curve on the smooth surface is a geodesic curve (i.e. $k_g(s) = 0$) then $\theta = \pi/2$, $k_n(s) = k(s)$ and equation of the rectifying curve is given by
\begin{equation}\begin{split}
\gamma(s)& = \lambda(s)(u^{'}\phi_u + v^{'}\phi_v) + \frac{\mu(s)}{\sqrt{EG-F^2}}\\&(Eu^{'}\phi_v + F(v^{'}\phi_v - u^{'}\phi_u) -G v^{'}\phi_u). 
\end{split}\end{equation} 

\n Also, if the rectifying curve on the smooth surface is a asymptotic line (i.e. $k_n(s) = 0$) then $\theta = 0$, $k_g(s) = k(s)$ and equation of the rectifying curve is given by
\begin{equation}
\gamma(s) = \lambda(s)(u^{'}\phi_u + v^{'}\phi_v) + \frac{\mu(s)}{\sqrt{EG-F^2}}(\phi_u\times\phi_v). 
\end{equation} 

\begin {theorem} Let $f: S\rightarrow \overline{S}$ be an isometry, where $S$ and $\overline{S}$ are smooth surfaces and $\gamma(s)$ be a rectifying curve on $S$ with $k_n \neq 0$. Then $\overline{\gamma} = fo\gamma$ is a rectifying curve on $\overline{S}$ if any one of the following conditions holds:

(i) $\overline\gamma$ is geodesic curve on $\overline S$ and $\overline{\gamma} (s) = f_*(\gamma(s))$,

(ii) $\overline{\gamma}$ is asymptotic curve on $\overline{S}$ and $\overline{\gamma}(s) + \dfrac{\mu(s)k_n(s)}{k(s)}\overline{P}  = f_*(\gamma(s))$,

(iii) $\overline{\gamma}$ is neither geodesic nor asymptotic curve on $\overline{S}$ and $\overline{\gamma} (s) = f_*(\gamma(s))$.

\end{theorem}

\n \textit{Proof.} Let $f: S\rightarrow \overline{S}$ be an isometry, where $S$ and $\overline{S}$ are smooth surfaces and $\gamma(s)$ be a rectifying curve on $S$ such that $k_n \neq 0$.
 
\n Suppose (i) holds. Then, $\overline k_g(s) = 0$ and $\overline{\gamma}(s) = f_*(\gamma(s))$, which implies
\begin{equation}\begin{split}
\overline{\gamma}(s) =& \lambda(s)(u^{'}f_*{\phi_u} + v^{'}f_*{\phi_v}) + \frac{\mu(s)}{\sqrt{EG-F^2}}(Eu^{'}f_*{\phi_v} + F(v^{'}f_*{\phi_v} - u^{'}f_*{\phi_u})\\& - Gv^{'}f_*{\phi_u}).
\end{split}\end{equation}
Thus from (2.13) and (3.6), we get
\begin{equation}
\overline{\gamma}(s) = \overline\lambda(s)(u^{'}\overline{\phi_u} + v^{'}\overline{\phi_v}) + \frac{\overline\mu(s)}{\sqrt{\overline E\overline G-\overline F^2}}(\overline Eu^{'}\overline{\phi_v} + \overline F(v^{'}\overline{\phi_v} - u^{'}\overline{\phi_u}) - \overline Gv^{'}\overline{\phi_u}),
\end{equation}
 where $\overline{\lambda}(s) = \lambda(s)$ and $\overline{\mu}(s) = \mu(s)$. This is equation of rectifying curve on $\overline S$, which is geodesic on the surface.

\n Suppose (ii) holds. Then, $\overline k_n(s) = 0$ and $\overline{\gamma}(s) = f_*(\gamma(s)) - \dfrac{\mu(s)k_n(s)}{k(s)}\overline{P}$, which implies
\begin{equation}\begin{split}
\overline{\gamma}(s) &= \lambda(s)(u^{'}f_*{\phi_u} + v^{'}f_*{\phi_v}) +\dfrac{\mu(s)k_g(s)}{k(s)}f_*{U} - \dfrac{\mu(s)k_n(s)}{k(s)}\overline{P} \\&+ \dfrac{\mu(s)k_n(s)}{k(s)\sqrt{EG-F^2}}(Eu^{'}f_*{\phi_v} + F(v^{'}f_*{\phi_v} - u^{'}f_*{\phi_u}) - Gv^{'}f_*{\phi_u}).
\end{split}\end{equation}
Thus from (2.13) and (3.8), we obtain
\begin{equation}\begin{split}
\overline{\gamma}(s) &= \lambda(s)(u^{'}\overline{\phi_u} + v^{'}\overline{\phi_v}) +\dfrac{\mu(s)k_g(s)}{k(s)}\overline{U} - \dfrac{\mu(s)k_n(s)}{k(s)}\overline{P} \\&+ \dfrac{\mu(s)k_n(s)}{k(s)\sqrt{\overline E\overline G-\overline F^2}}(\overline Eu^{'}\overline{\phi_v} + \overline F(v^{'}\overline{\phi_v} - u^{'}\overline{\phi_u}) - \overline Gv^{'}\overline{\phi_u}),
\end{split}\end{equation}
which implies
\begin{equation}
\overline{\gamma}(s) = \overline\lambda(s)(u^{'}\overline{\phi_u} + v^{'}\overline{\phi_v}) +\overline\mu(s)\overline{U},
\end{equation}
where $\overline{\lambda}(s) = \lambda(s)$ and $\overline{\mu}(s) = \dfrac{\mu(s)k_g(s)}{k(s)}$.
This is equation of rectifying curve on $\overline S$, which is asymptotic on the surface.

\n Now, suppose (iii) holds. Then,  $\overline k_g(s) \neq 0$, $\overline k_n(s) \neq 0$ and $\overline{\gamma}(s) = f_*(\gamma(s))$, which implies
\begin{equation}\begin{split}
\overline{\gamma}(s) = &\lambda(s)(u^{'}f_*{\phi_u} + v^{'}f_*{\phi_v}) +\dfrac{\mu(s)k_g(s)}{k(s)}f_*{U} \\& + \dfrac{\mu(s)k_n(s)}{k(s)\sqrt{EG-F^2}}(Eu^{'}f_*{\phi_v} + F(v^{'}f_*{\phi_v} - u^{'}f_*{\phi_u}) - Gv^{'}f_*{\phi_u}).
\end{split}\end{equation}
Thus from (2.13) and (3.11), we get
\begin{equation}\begin{split}
\overline{\gamma}(s) =& \lambda(s)(u^{'}\overline{\phi_u} + v^{'}\overline{\phi_v}) +\dfrac{\mu(s)k_g(s)}{k(s)}\overline{U} \\& + \dfrac{\mu(s)k_n(s)}{k(s)\sqrt{\overline E\overline G-\overline F^2}}(\overline Eu^{'}\overline{\phi_v} + \overline F(v^{'}\overline{\phi_v} - u^{'}\overline{\phi_u}) - \overline Gv^{'}\overline{\phi_u}),
\end{split}\end{equation}
which implies
\begin{equation}\begin{split}
\overline{\gamma}(s) = &\overline\lambda(s)(u^{'}\overline{\phi_u} + v^{'}\overline{\phi_v}) +\dfrac{\overline\mu(s)\overline{k_g}(s)}{\overline{k}(s)}\overline{U} \\&+ \dfrac{\overline\mu(s)\overline{k_n}(s)}{\overline{k}(s)\sqrt{\overline E\overline G-\overline F^2}}(\overline Eu^{'}\overline{\phi_v} + \overline F(v^{'}\overline{\phi_v} - u^{'}\overline{\phi_u}) - \overline Gv^{'}\overline{\phi_u}),
\end{split}\end{equation} 
where $\overline{\lambda}(s) = \lambda(s)$, $\dfrac{\overline\mu(s)\overline{k_g}(s)}{\overline{k}(s)} = \dfrac{\mu(s)k_g(s)}{k(s)}$ and $\dfrac{\overline\mu(s)\overline{k_n}(s)}{\overline{k}(s)} = \dfrac{\mu(s)k_n(s)}{k(s)}$.

 \n This is equation of rectifying curve on $\overline S$, which is neither geodesic nor asymptotic on the surface.

\begin {theorem} Let $f: S\rightarrow \overline{S}$ be an isometry, where $S$ and $\overline{S}$ are smooth surfaces and $\gamma(s)$ be a rectifying curve on $S$ with $k_n = 0$. Then $\overline{\gamma} = fo\gamma$ is a rectifying curve on $\overline{S}$ if any one of the following conditions holds:

(i) $\overline{\gamma}$ is asymptotic curve on $\overline{S}$ and $\overline{\gamma}(s) = F_*(\gamma(s))$,

(ii) $\overline{\gamma}$ is not asymptotic curve on $\overline{S}$ and $\overline{\gamma} (s) - \mu(s)\overline{P}  = F_*(\gamma(s))$.

\end{theorem}

\n \textit{Proof.} We can easily prove by using Theorem 3.1.

\begin {theorem} Let $f: S\rightarrow \overline{S}$ be an isometry. If $\gamma$ and $\overline{\gamma}$ are rectifying curves on $S$ and $\overline{S}$ respectively, with $k_n \neq 0$ then, we have follwing:

(i) if $\overline \gamma$ is geodesic curve on $\overline S$ then $\overline{\gamma}(s).\overline{T}(s) = \gamma(s).T(s)$,

(ii) if $\overline{\gamma}$ is asymptotic curve on $\overline{S}$ then  $\overline{\gamma}(s).\overline{T}(s) - \gamma(s).T(s) = \dfrac{\mu(s)k_n(s)}{k(s)})\\(\sqrt{EG-F^2})(av^{'} - bu^{'})$,

(iii) if $\overline{\gamma}$ is neither geodesic nor asymptotic curve on $\overline{S}$ then $\overline{\gamma}(s).\overline{T}(s) = \gamma(s).T(s)$,

where $T(s) = a\phi_u + b\phi_v$ is any tangent vector to the surface $S$ at point $\gamma(s)$.
\end{theorem}

\n \textit{Proof.} Let $f: S\rightarrow \overline{S}$ be an isometry and $\gamma$ and $\overline{\gamma}$ be rectifying curves on $S$ and $\overline{S}$ respectively, with $k_n \neq 0$. Then,
\begin{equation}
\overline{\gamma}(s).\overline{T}(s) - \gamma(s).T(s) = a(\overline{\gamma}(s).\overline\phi_u - \gamma(s).\phi_u) + b(\overline{\gamma}(s).\overline\phi_v  - \gamma(s).\phi_v).
\end{equation}
Now, from (3.3), we get
\begin{equation}
\gamma(s).\phi_u = \lambda(s)(u^{'}E + v^{'}F) + \dfrac{\mu(s)k_n(s)}{k(s)\sqrt{EG-F^2}}(F^{2} - GE)v^{'}.
\end{equation}
Similarly, we obtain
\begin{equation}
\gamma(s).\phi_v = \lambda(s)(u^{'}F + v^{'}G) - \dfrac{\mu(s)k_n(s)}{k(s)\sqrt{EG-F^2}}(F^{2} - GE)u^{'}.
\end{equation}
For (i), suppose $\overline \gamma$ is geodesic curve on $\overline S$. Then from (2.13) and (3.4), we get
\begin{equation}
\gamma(s).\phi_u = \lambda(s)(u^{'}E + v^{'}F) + \dfrac{\mu(s)}{\sqrt{EG-F^2}}(s)(F^{2} - GE)v^{'},
\end{equation}
and
\begin{equation}
\overline\gamma(s).\overline\phi_u = \lambda(s)(u^{'}E + v^{'}F) + \dfrac{\overline\mu(s)}{\sqrt{EG-F^2}}(F^{2} - GE)v^{'}.
\end{equation}
Thus from (3.17) and (3.18), we obtain
\begin{equation}
\overline\gamma(s).\overline\phi_u - \gamma(s).\phi_u = (\mu(s) - \overline\mu(s))(\sqrt{EG-F^2})v^{'}.
\end{equation}
Similarly, we get
\begin{equation}
\overline\gamma(s).\overline\phi_v - \gamma(s).\phi_v = (\mu(s) - \overline\mu(s))(\sqrt{EG-F^2})u^{'}.
\end{equation}
Thus from (3.14), (3.19) and (3.20), we get
\begin{equation}
\overline{\gamma}(s).\overline{T}(s) - \gamma(s).T(s) = (\mu(s) - \overline\mu(s))(\sqrt{EG-F^2})(av^{'} + bu^{'}).
\end{equation}
Since $\gamma$ and $\overline{\gamma}$ are rectifying curves on $S$ and $\overline{S}$, respectively therefore $\overline\mu(s) = \mu(s)$. Hence, $\overline{\gamma}(s).\overline{T}(s) = \gamma(s).T(s)$.

\n For (ii), suppose $\overline{\gamma}$ is asymptotic curve on $\overline{S}$. Then from (2.13) and (3.5), we get
\begin{equation}
\overline\gamma(s).\overline\phi_u = \lambda(s)(u^{'}E + v^{'}F),
\end{equation}
and
\begin{equation}
\overline\gamma(s).\overline\phi_v = \lambda(s)(u^{'}F + v^{'}G).
\end{equation}
Thus from (3.15) and (3.22), we obtain
\begin{equation}
\overline\gamma(s).\overline\phi_u - \gamma(s).\phi_u = - \dfrac{\mu(s)k_n(s)}{k(s)\sqrt{EG-F^2}}(F^{2} - GE)v^{'}.
\end{equation}
Also from (3.16) and (3.23), we get
\begin{equation}
\overline\gamma(s).\overline\phi_v - \gamma(s).\phi_v =    \dfrac{\mu(s)k_n(s)}{k(s)\sqrt{EG-F^2}}(F^{2} - GE)u^{'}.
\end{equation}
Thus from (3.14), (3.24) and (3.25), we obtain
\begin{equation}
\overline{\gamma}(s).\overline{T}(s) - \gamma(s).T(s) =  \dfrac{\mu(s)k_n(s)}{k(s)})(\sqrt{EG-F^2})(av^{'} - bu^{'}).
\end{equation}
  
\n Now for (iii), suppose $\overline{\gamma}$  is neither geodesic nor asymptotic curve on $\overline{S}$. Then from (2.13) and (3.3), we get
\begin{equation}
\overline\gamma(s).\overline\phi_u = \lambda(u^{'}E + v^{'}F) + \dfrac{\overline\mu(s)\overline k_n(s)}{\overline k(s)\sqrt{EG-F^2}}(F^{2} - GE)v^{'},
\end{equation}
and
\begin{equation}
\overline\gamma(s).\overline\phi_v = \lambda(u^{'}F + v^{'}G) - \dfrac{\overline\mu(s)\overline k_n(s)}{\overline k(s)\sqrt{EG-F^2}}(F^{2} - GE)u^{'}.
\end{equation}
Thus from (3.15) and (3.27), we obtain
\begin{equation}
\overline\gamma(s).\overline\phi_u - \gamma(s).\phi_u = (\dfrac{\mu(s)k_n(s)}{k(s)} - \dfrac{\overline\mu(s)\overline k_n(s)}{\overline k(s)})\sqrt{EG-F^2})v^{'}.
\end{equation}
Similarly, from (3.16) and (3.28), we get
\begin{equation}
\overline\gamma(s).\overline\phi_v - \gamma(s).\phi_v = (\dfrac{\overline\mu(s)\overline k_n(s)}{\overline k(s)} - \dfrac{\mu(s)k_n(s)}{k(s)})(\sqrt{EG-F^2})u^{'}.
\end{equation}
Thus, from (3.14), (3.29) and (3.30), we obtain
\begin{equation}
\overline{\gamma}(s).\overline{T}(s) - \gamma(s).T(s) = (\dfrac{\overline\mu(s)\overline k_n(s)}{\overline k(s)} - \dfrac{\mu(s)k_n(s)}{k(s)})(\sqrt{EG-F^2})(bu^{'} - av^{'}).
\end{equation}
Since $\gamma$ and $\overline{\gamma}$ are rectifying curves on $S$ and $\overline{S}$, respectively therefore $\overline{\gamma}(s).\overline{T}(s) = \gamma(s).T(s)$.

 \begin {theorem} Let $f: S\rightarrow \overline{S}$ be an isometry. If $\gamma$ and $\overline{\gamma}$ are rectifying curves on $S$ and $\overline{S}$ respectively, with $k_n = 0$ then, we have follwing:
 
(i) if $\overline{\gamma}$ is asymptotic curve on $\overline{S}$ then $\overline{\gamma}(s).\overline{T}(s) = \gamma(s).T(s)$,

(ii) if $\overline{\gamma}$ is not asymptotic curve on $\overline{S}$ then  $\overline{\gamma}(s).\overline{T}(s) - \gamma(s).T(s) = \dfrac{\overline\mu(s)\overline k_n(s)}{\overline k(s)}\\(\sqrt{EG-F^2})(bu^{'} - av^{'})$.
\end{theorem}

\n \textit{Proof.} We can easily prove by using Theorem 3.3.

\begin {theorem}  Let $f: S\rightarrow \overline{S}$ be an isometry. If $\gamma$ and $\overline{\gamma}$ are rectifying curves on $S$ and $\overline{S}$ respectively, with $k_n \neq 0$ then,  we have follwing:

(i) if $\overline\gamma$ is geodesic curve on $\overline S$ then $\overline{\gamma}(s).\overline{P}(s) = \gamma(s).P(s)$,

(ii) if $\overline{\gamma}$ is asymptotic curve on $\overline{S}$ then  $\overline{\gamma}(s).\overline{P}(s) - \gamma(s).P(s) = ((aE + bF)u^{'}  - (aF + bG)v^{'})(\sqrt{EG-F^2}) \dfrac{\mu(s)k_n(s)}{k(s)}$,

(iii) if $\overline{\gamma}$ is neither geodesic nor asymptotic curve on $\overline{S}$ then $\overline{\gamma}(s).\overline{P}(s) = \gamma(s).P(s)$, 

where $P = U\times T$ and $T = a\phi_u + b\phi_v$ be any tangent vector to the surface $S$ at point $\gamma(s)$.
\end{theorem}
\n \textit{Proof.} Let $f: S\rightarrow \overline{S}$ be an isometry. If $\gamma$ and $\overline{\gamma}$ are rectifying curves on $S$ and $\overline{S}$ respectively, with $k_n \neq 0$ then, by using (2.11), we obtain
\begin{equation}
\gamma(s).P(s) = \dfrac{(aE + bF)}{\sqrt{EG-F^2}}\gamma(s).\phi_{v} - \dfrac{(aF + bG)}{\sqrt{EG-F^2}}\gamma(s).\phi_{u}.
\end{equation}
Similarly by using (2.13), we get
\begin{equation}
\overline\gamma(s).\overline P(s) = \dfrac{(aE + bF)}{\sqrt{EG-F^2}}\overline\gamma(s).\overline\phi_{v} - \dfrac{(aF + bG)}{\sqrt{EG-F^2}}\overline\gamma(s).\overline\phi_{u}.
 \end{equation}
 Thus from (3.32) and (3.33), we obtain
 \begin{equation}\begin{split}
 \overline\gamma(s).\overline P(s) -\gamma(s).P(s) =& \dfrac{(aE + bF)}{\sqrt{EG-F^2}}(\overline\gamma(s).\overline\phi_{v} - \gamma(s).\phi_{v}) \\&-\dfrac{(aF + bG)}{\sqrt{EG-F^2}}(\overline\gamma(s).\overline\phi_{u} - \gamma(s).\phi_{u}).
 \end{split}\end{equation}
 
 \n For (i), suppose $\overline\gamma$ is geodesic curve on $\overline S$. Then from (3.19), (3.20) and (3.34), we get
 \begin{equation}\begin{split}
 \overline\gamma(s).\overline P(s) -\gamma(s).P(s) = &((aE + bF)u^{'}  + (aF + bG)v^{'})(\sqrt{EG-F^2})\\&(\overline\mu(s) - \mu(s)).
 \end{split}\end{equation}
 Since $\gamma$ and $\overline{\gamma}$ are rectifying curves on $S$ and $\overline{S}$, respectively therefore $\overline\mu(s) = \mu(s)$. Hence, $\overline\gamma(s).\overline P(s) = \gamma(s).P(s)$.
 
 \n For (ii), suppose $\overline{\gamma}$ is asymptotic curve on $\overline{S}$. Then from (3.24), (3.25) and (3.34), we obtain
 \begin{equation}\begin{split}
 \overline\gamma(s).\overline P(s) -\gamma(s).P(s) = &((aE + bF)u^{'}  - (aF + bG)v^{'})(\sqrt{EG-F^2})\\& \dfrac{\mu(s)k_n(s)}{k(s)}.
 \end{split}\end{equation}
 
\n Now for (iii), suppose $\overline{\gamma}$  is neither geodesic nor asymptotic curve on $\overline{S}$. Then from (3.29), (3.30) and (3.34), we get
 \begin{equation}\begin{split}
 \overline\gamma(s).\overline P(s) -\gamma(s).P(s) = &((aE + bF)u^{'}  + (aF + bG)v^{'})(\sqrt{EG-F^2})\\&(\dfrac{\overline\mu(s)\overline k_n(s)}{\overline k(s)} - \dfrac{\mu(s)k_n(s)}{k(s)}).
 \end{split}\end{equation}
 Since $\gamma$ and $\overline{\gamma}$ are rectifying curves on $S$ and $\overline{S}$, respectively therefore $\dfrac{\overline\mu(s)\overline k_n(s)}{\overline k(s)} = \dfrac{\mu(s)k_n(s)}{k(s)}$. Hence, $\overline{\gamma}(s).\overline{P}(s) = \gamma(s).P(s)$.
 
\begin {theorem}  Let $f: S\rightarrow \overline{S}$ be an isometry. If $\gamma$ and $\overline{\gamma}$ are rectifying curves on $S$ and $\overline{S}$ respectively, with $k_n = 0$ then, we have follwing:

(i) if $\overline{\gamma}$ is asymptotic curve on $\overline{S}$ then  $\overline{\gamma}(s).\overline{P}(s) = \gamma(s).P(s)$,

(ii) if $\overline{\gamma}$ is not asymptotic curve on $\overline{S}$ then $\overline{\gamma}(s).\overline{P}(s) - \gamma(s).P(s) = ((aE + bF)u^{'}  + (aF + bG)v^{'})(\sqrt{EG-F^2})(\dfrac{\overline\mu(s)\overline k_n(s)}{\overline k(s)})$, 
where $P = U\times T$ and $T = a\phi_u + b\phi_v$ be any tangent vector to the surface $S$ at point $\gamma(s)$.
\end{theorem}
\n \textit{Proof.} We can easily prove by using Theorem 3.5.

\begin {theorem}  Let $f: S\rightarrow \overline{S}$ be an isometry. If $\gamma$ and $\overline{\gamma}$ are rectifying curves on $S$ and $\overline{S}$ respectively then, we have  $\overline{\gamma}(s).\overline{U}(s) = \gamma(s).U(s)$.
\end{theorem}

\n \textit{Proof.}  Let $f: S\rightarrow \overline{S}$ be an isometry. Then from (3.3), we have \begin{equation}
\overline\gamma(s).\overline U(s) = \overline\gamma(s).(\dfrac{\overline\phi_{u}\times\overline\phi_{v}}{\sqrt{EG-F^2}}) = \dfrac{\overline\mu(s)\overline k_g(s)}{k(s)},
\end{equation}
and
\begin{equation}
\gamma(s).U(s) = \gamma(s).(\dfrac{\phi_{u}\times\phi_{v}}{\sqrt{EG-F^2}}) = \dfrac{\mu(s)k_g(s)}{k(s)}.
\end{equation}
Thus if $k_g = 0$ then $\overline{\gamma}(s).\overline{U}(s) = \gamma(s).U(s)$.

\n Also, if $k_g \neq 0$ then $\dfrac{\overline\mu(s)\overline k_g(s)}{\overline k(s)} = \dfrac{\mu(s)k_g(s)}{k(s)}$ and hance $\overline{\gamma}(s).\overline{U}(s) = \gamma(s).U(s)$.

\n Note:  Let $f: S\rightarrow \overline{S}$ be an isometry. If $\gamma$ and $\overline{\gamma}$ are rectifying curves on $S$ and $\overline{S}$, respectively then the components of $\gamma(s)$ along $T = a\phi_u + b\phi_v$, $U$ and $P = U\times T$ are invariant under isometry if any one of the following holds:

(i) both $\gamma$ and $\overline{\gamma}$ are asymptotic on $S$ and $\overline{S}$ respectively,

(ii) both $\gamma$ and $\overline{\gamma}$ are geodesic on $S$ and $\overline{S}$ respectively,

(iii) neither $\gamma$ nor $\overline{\gamma}$ are geodesic and asymptotic on $S$ and $\overline{S}$ respectively.

\noindent\author{Akhilesh Yadav, Buddhadev Pal}\\
\date{Department of Mathematics, Institute of Science, \\Banaras Hindu University, Varanasi-221005, India}\\
\maketitle {\noindent E-mails: akhilesh$\_$mathau@rediffmail.com, pal.buddha@gmail.com} 


\begin{thebibliography}{}
      

\bibitem[1]{} Chen, B. Y., \textit{When does the position vector of a space curve always lie in its rectifying plane?}, Amer.
Math. Monthly, 110 (2003), 147-152.

\bibitem[2]{} Chen, B. Y., \textit{Rectifying curves and geodesics on a cone in the Euclidean 3-space}, Tamkang J. Math., 48 (2017), 209-214.

\bibitem[3]{} Chen, B. Y., Dillen, F., \textit{Rectifying curve as centrode and extremal curve}, Bull. Inst. Math. Acad.
Sinica, 33(2) (2005), 77-90.

\bibitem[4]{} Camci, C., Kula, L. and Ilarslan, K., \textit{Characterizations of the position vector of a surface curve in Euclidean 3-space}, An. St. Univ. Ovidius Constanta, 19, no. 3, (2011), 59-70.

\bibitem[5]{} Deshmukh, S., Chen, B. Y. and Alshammari, S. H., \textit{On a rectifying curves in Euclidean 3-space}, Turk. J.
Math., 42 (2018), 609-620.

\bibitem[6]{} do Carmo, M. P., \textit{Differential geometry of curves and surfaces}, Prentice-Hall, Inc, New Jersey, (1976). 

\bibitem[7]{} Ilarslan, K., Nesovic, E. \textit{Some characterizations of rectifying curves in the Euclidean space $E^4$}, Turk. J. Math., 32 (2008), 21-30.

\bibitem[8]{} Kim D. S., Chung, H. S. and Cho, K. H., \textit{Space curves satisfying $\tau/\kappa= as + b$}, Honam Math. J., 15 (1993), 5-9.

\bibitem[9]{} Pressley, A., \textit{Elementary differential geometry}, Springer-Verlag, 2001.

\bibitem[10]{} Shaikh, A. A., Ghosh, P. R. \textit{Rectifying curves on a smooth surface immersed in the Euclidean space},
Indian J. Pure Appl. Math., 50(4) (2019), 883-890.


\end{thebibliography}
\end{document}